\documentclass[12pt]{article}

\usepackage{amssymb}

\newtheorem{lemma}{Lemma}[section]
\newtheorem{theorem}{Theorem}[section]

\newtheorem{proposition}{Proposition}[section]
\newtheorem{corollary}{Corollary}[section]

\begin{document}



\title{Sharp Nash inequalities on manifolds with boundary in the presence of symmetries}


\author{{Athanase Cotsiolis and Nikos Labropoulos}\\
\small{\textit{Department of Mathematics University of Patras}},\\
\small{\textit{Patras 26110, Greece}}\\
\scriptsize{e-mails: cotsioli@math.upatras.gr and nal@upatras.gr}}
 \maketitle
\noindent\textbf{Abstract }\\
In this paper we establish the best constant $\widetilde
A_{opt}(\overline{M})$ for the Trace Nash inequality on a
$n-$dimensional compact Riemannian manifold in the presence of
symmetries, which is an improvement over the classical case due to
the symmetries which arise and reflect the geometry of manifold.
This is particularly true when the data of the problem is
invariant under the action of an arbitrary compact subgroup $G$ of
the isometry group $Is(M,g)$, where all the orbits have infinite
cardinal.\\





\section{Introduction}
$\;\;\;\;$ We say that the Nash inequality (\ref{E1.1}) is valid
if there exists a constant $A>0$ such that for all $u \in
C_0^\infty \left( \mathbb{R}^n \right)$, $n\geq 2$
\begin{equation}\label{E1.1}
 \left( {\int_{\mathbb{R}^n} {u^2 dx } }
\right)^{1 + \frac{2} {n}} \leqslant A {\int_{\mathbb{R}^n}{\left|
{\nabla u} \right|^2 } dx   } \left( {\int_{\mathbb{R}^n} {\left|
u \right|dx } } \right)^{\frac{4} {n}}
\end{equation}
Such an inequality first appeared in the celebrated paper of Nash
\cite{Nas}, where he discussed the H\"older regularity of
solutions of divergence form in uniformly elliptic equations. It
is a particular case of the Gagliardo-Nirenberg type inequalities
$||u||_r  \leqslant C ||\nabla u||_q^a \,||u||_s^{1-a}$ and it is
well known that the Nash inequality (\ref{E1.1}) and the Euclidian
type Sobolev inequality are equivalent in the sense that if one of
them is valid, the other one is also valid (i.e. see
\cite{Bak-Cou-Led-Sal}). It is, also, well known  that with this
procedure of passing from the one type of inequalities to the
other, is impossible to compare  the best constants, since the inequalities under use are not optimal. \\
As far as the optimal version of Nash inequality (\ref{E1.1}) is
concerned, the best constant $A_0(n)$, that is
\[
A_0 \left( n \right)^{ - 1}  = \inf \left\{ {\left.
{\frac{{\int_{\mathbb{R}^n } {\left| {\nabla u} \right|^2 dx}
\left( {\int_{\mathbb{R}^n } {\left| u \right|dx} }
\right)^{\frac{4} {n}} }} {{\left( {\int_{\mathbb{R}^n } {u^2 dx}
} \right)^{1 + \frac{2} {n}} }}} \right|u \in C_0^\infty  \left(
{\mathbb{R}^n } \right),u\not \equiv 0} \right\},
\]
has been computed by Carlen and Loss in \cite{Car-Los}, together
with the characterization of the extremals for the corresponding
optimal inequality, as
\[
A_0 \left( n \right) = \frac{{\left( {n + 2} \right)^{\frac{{n +
2}} {n}} }} {{2^{\frac{2} {n}} n\lambda _1^N \left| {\mathcal{B}^n
} \right|^{\frac{2} {n}} }},
\]
where $\left| {\mathcal{B}^n } \right|$ denotes the euclidian
volume of the unit ball $\mathcal{B}^n $ in $\mathbb{R}^n$ and
$\lambda _1^N$ is the first Neumann eigenvalue for the Laplacian
for radial functions in the  unit ball $\mathcal{B}^n
$.\\
For an example of application of the Nash inequality with the best
constant, we refer to Kato \cite{Kat} and for a geometric proof
with an asymptotically sharp constant, we refer to Beckner
\cite{Bec}.\\

For compact Riemannian manifolds, the Nash inequality still holds
with an additional $L^1-$term and that is why we will refer this
as the $L^1-$Nash inequality.\\
Given $(M, g)$ a smooth compact Riemannian $n-$manifold, $n\geq
2$, we get here the existence of real constants $A$ and $B$ such
that for any $u \in C^\infty(M)$,
\begin{eqnarray}\label{E1.2}
 \left( {\int_M {u^2 dV_g } }
\right)^{1 + \frac{2} {n}} \leqslant A  {\int_M {\left| {\nabla u}
\right|_g^2 } dV_g   } \left( {\int_M {\left| u \right|dV_g } }
\right)^{\frac{4} {n}}+B\left( {\int_M {\left| u \right|dV_g } }
\right)^{2+\frac{4} {n}}
\end{eqnarray}
The best constant for this inequality is defined as
\[
A^1_{opt}(M) = \inf \left\{ { A > 0:\exists \, B
> 0\,\,\,\mathrm{s.t.}\,\,(\ref{E1.2}) \,\,\,
\mathrm{is\,\,true}\,\forall \,u \in C^\infty \left( M \right)}
\right\}
\]
This inequality has been  studied completely by Druet, Hebey and
Vaugon. They proved in \cite {Dru-Heb-Vau} that
$A^1_{opt}(M)=A_0(n)$, and (\ref{E1.2}) with its optimal constant
$A=A_0(n)$ is sometimes valid and sometimes not, depending on the
geometry of $M$.\\
Humbert in \cite{Hum1} studied the following $L^2-$Nash inequality
\begin{equation}\label{E1.3}
\left( {\int_M {u^2 dV_g } } \right)^{1 + \frac{2} {n}} \leqslant
\left( {A\int_M {\left| {\nabla u} \right|_g^2 } dV_g  + B\int_M
{u^2 dV_g } } \right)\left( {\int_M {\left| u \right|dV_g } }
\right)^{\frac{4} {n}},
\end{equation}
for all $u \in C^\infty  \left( M \right)$, of which the best
constant is  defined as
\[
A^2_{opt}(M)  = \inf \left\{ {A > 0:\exists \,B >
0\,\,\,\mathrm{s.t.}\,\,\,(\ref{E1.3})
\,\,\,\mathrm{is\,\,true}\,\forall \,u \in C^\infty  \left( M
\right)} \right\}
 \]
Contrary to the sharp $L^1-$Nash inequality, in this case, he
proved that $B$ always exists and $A^2_{opt}(M)=A_0(n)$.\\

We denote $\, \mathbb{R}^n_+=\mathbb{R}^{n-1}\times[0,
+\infty)\,\, \mathrm{and}\,\,  \partial
\mathbb{R}^n_+=\mathbb{R}^{n-1}\times\{0\} $. The trace Nash
inequality  states that a constant $\tilde{A}>0$ exists such that
for all $u \in C_0^\infty (\mathbb{R}^n_+)$, $n\geq 2$ with
$\nabla u \in L^2(\mathbb{R}^n)$ and $u|_{\partial
\mathbb{R}^n_+}\in L^1(\partial\mathbb{R}^n_+ )\cap
L^2(\partial\mathbb{R}^n_+ )$
\begin{equation}\label{E1.4}
 \left( {\int_{\partial\mathbb{R}^n_+} {u^2 ds } }
\right)^ \frac{n} {n-1} \leqslant \widetilde{A}
{\int_{\mathbb{R}^n_+} {\left| {\nabla u} \right|^2 } dx   }
\left( {\int_{\partial\mathbb{R}^n_+} {\left| u \right|ds } }
\right)^{\frac{2} {n-1}},
\end{equation}
where $ds$ is the standard volume element on
$\mathbb{R}^{n-1} $ and the trace of $u$ on $\partial \mathbb{R}^n_+$ is also denoted by  $u$.\\
Let $\widetilde{A}_0(n)$ be the best constant in Nash inequality
 (\ref{E1.4}).  That is
\[
\widetilde A_0 \left( n \right)^{ - 1}  = \inf \left\{ {\left.
{\frac{{\int_{ \mathbb{R}_ + ^n } {\left| {\nabla u} \right|^2 dx}
\left( {\int_{\partial \mathbb{R}_ + ^n } {\left| u \right|ds} }
\right)^{\frac{2} {n-1}} }} {\left({\int_{\partial \mathbb{R}_ +
^n } {u^2 ds} }\right)^{\frac{n}{n-1}}}} \right|u \in C_0^\infty
\left( {\mathbb{R}_ + ^n } \right),u\not  \equiv 0} \right\}
\]
The computation problem of the exact value of $\widetilde A_0
\left( n \right)$ still remains open.\\

For compact Riemannian manifolds with boundary, Humbert, also,
studied in \cite{Hum2} the trace  Nash inequality.\\
On smooth compact  $n-$dimensional, $n\geq 2$, Riemannian
manifolds with boundary, for all $u \in C^\infty(M)$, consider the
following trace Nash inequality
\begin{equation}\label{E1.5}
 \left( {\int_{\partial M} {u^2 dS_g } }
\right)^{\frac{n} {{n - 1}}}  \leqslant \left( {\tilde A\int_M
{\left| {\nabla u} \right|_g^2 } dV_g  + \tilde B\int_{\partial M}
{u^2 dS_g } } \right)\left( {\int_{\partial M} {\left| u
\right|dS_g } } \right)^{\frac{2} {{n - 1}}}
\end{equation}
The best constant for the above inequality is defined as
\[
\widetilde A_{opt}(M) = \inf \left\{ {\tilde A > 0:\exists
\,\tilde B
> 0\,\,\,\mathrm{s.t.}\,\,(\ref{E1.5}) \,\,\,
\mathrm{is\,\,true}\,\forall \,u \in C^\infty \left( M \right)}
\right\}
\]
It was proved in \cite {Hum2} that $\widetilde
A_{opt}(M)=\widetilde A_0(n)$, and (\ref{E1.5}) with its optimal
constant
$\widetilde A=\widetilde A_0(n)$ is always valid.\\

In this paper we prove that, when the functions are invariant
under an isometry group, all orbits of which are of infinite
cardinal, the Nash inequalities can be improved, in the sense that
we can get a higher critical
exponent.\\ More precisely we establish:\\
\textbf{(A)}$\,\,$ The best constant for the \textbf{Nash
inequality} on compact Riemannian manifolds with boundary,
invariant under the action of an arbitrary compact subgroup $G$ of
the isometry group $Is(M,g)$, where all the orbits have infinite
cardinal, and\\
\textbf{(B)}$\,\,$ The best constant for the \textbf{Trace Nash
inequality} on compact Riemannian manifolds with boundary,
invariant under the action of an arbitrary compact subgroup $G$ of
the isometry group
$Is(M,g)$, where all the orbits have infinite cardinal.  \\
These best constants are improvements over the classical cases due
to the symmetries which arise and reflect the geometry of the
manifold.\\

\section{Results and Examples }
\subsection{Results}
\begin{theorem}\label{T2.1}
Let $(M,g)$ be a smooth, compact $n-$dimensional Riemannian
manifold, $n\geq 3$, with boundary, $G-$invariant under the action
of a subgroup $\,G$ of the isometry group $Is(M,g)$. Let $k$
denotes the minimum orbit dimension of $G$ and $V$ denotes the
minimum of the volume of the $k-$dimension orbits. Then for any
$\varepsilon>0$ there exists a constant $B_\varepsilon$ such that
and for all $ u \in  H^2_{1,G}( M ) $ the following inequality
\begin{eqnarray*}
\left( {\int_M {u^2 dV_g } } \right)^{\frac{{n - k + 2}} {{n -
k}}} & \leqslant &  \left(\left( {A_G  + \varepsilon }
\right)^{\frac{{n - k}} {{n - k + 2}}} {\int_M {\left| {\nabla u}
\right|_g^2 } dV_g  + B_\varepsilon \int_M {u^2 dV_g } } \right)\\
&&\times\left( {\int_M {\left| u \right|dV_g } } \right)^{\frac{4}
{{n - k}}}
\end{eqnarray*}
holds, where $A_G  = \frac{{A_0( {n - k} )}}
{{V^{\frac{2} {{n - k }}} }}$.\\
Moreover the constants $ A_G $
is the best constant for this inequality.\\
\end{theorem}

\begin{theorem}\label{T2.2}
Let $(M,g)$ be a smooth, compact $n-$dimensional Riemannian
manifold, $n\geq 3$, with boundary, $G-$invariant under the action
of a subgroup $G$ of the isometry group $Is(M,g)$. Let $k$ denotes
the minimum orbit dimension of $G$ and $V$ denotes the minimum of
the volume of the $k-$dimension orbits. Then for any
$\varepsilon>0$ there exists a constant $\tilde B_\varepsilon$
such that and for all $ u \in  H^2_{1,G}( M ) $ the following
inequality
\begin{eqnarray*}
\left( {\int_{\partial M} {u^2 dS_g } } \right)^{\frac{{n - k}}
{{n - k - 1}}} & \leqslant &  \left(\left( {\tilde A_G  +
\varepsilon } \right)^{\frac{{n - k - 1}} {{n - k}}} {\int_M
{\left| {\nabla u} \right|^2 dV_g }
 +\tilde B_\varepsilon \int_{\partial M }{u^2 dS_g } } \right)\\
&&\times
 \left( {\int_{\partial M} {\left| u \right|} dS_g
} \right)^{\frac{2} {{n - k - 1}}}
\end{eqnarray*}
holds, where $\tilde A_G  = \frac{{\tilde A_0( {n - k} )}}
{{V^{\frac{1} {{n - k - 1}}} }}$.\\
Moreover the constants $ \tilde A_G $ is the best constant for
this inequality.\\
\end{theorem}

\begin{corollary}\label{T2.3} For any $\varepsilon>0$ there
exists a constant $C_\varepsilon$ such that and for all $u \in
H^2_{1,G} (T) $ the following inequality holds
$$
\left( {\int_{\partial T} {u^2 dS} } \right)^2  \le \left(
{\frac{{\tilde A_0 \left( 2 \right) + \varepsilon }}{{2\pi \left(
{l - r} \right)}}\int_T {\left| {\nabla u} \right|^2 dV}  +
C_\varepsilon \int_{\partial T} {u^2 dS} } \right)\left(
{\int_{\partial T} {\left| u \right|} dS} \right)^2
$$
Moreover the constant $\widetilde
A_{opt}(\overline{T})=\frac{{\tilde A_0 \left( 2 \right)}}{{2\pi
\left( {l - r} \right)}}$  is the best constant for this
inequality and verifies
\begin{eqnarray*}
\frac{3\sqrt{3}}{4 \pi^2(l-r)}\leq \widetilde
A_{opt}(\overline{T}) \leq \frac{2}{\pi^2(l-r)}
\end{eqnarray*}
\end{corollary}
\subsection{Examples}
 {\bf{Example 1.}}  Let
$T$ be the three dimensional solid torus
$$
T=\left\{(x,y,z)\in \mathbb{R}^3: ( \sqrt{x^2+y^2}-l)^2+z^2\leq
r^2,\, l>r>0\right\},
$$
with the metric induced by the $\mathbb{R}^3$ metric.
 Let $G=O(2)\times I$ be the group of rotations around axis
$z$. All $G-$orbits of $T$ are circles and the orbit of minimum
volume is the circle of radius $l-r$, and of length  $2\pi(l-r)$.
Then $T$ is a compact $3-$dimensional
 manifold with boundary, invariant under the action of
the subgroup $G$ of the isometry group $O(3)$.

In \cite{Cot-Lab-Tra} we found the best constant in inequality
(\ref{E1.5}) in the $3-$dimensional solid torus, which is
$G-$invariant under the action of a subgroup
$G=O\left( 2 \right)\times I$ of the isometry group $O(3)$.\\
{\bf{Example 2.}}  Let
$\mathbb{R}^n=\mathbb{R}^k\times\mathbb{R}^m $, $k\geq 2$, $m\geq
1$ and $\overline{\Omega}\subset(\mathbb{R}^k \backslash
\{0\})\times \mathbb{R}^m $. Let $G_{k,m}=O(k)\times Id_m$, the
subgroup of the isometry group $O(n)$ of the type $(x_1,
x_2)\longrightarrow (\sigma (x_1),x_2)$, $\sigma\in O(k)$, $x_1
\in \mathbb{R}^k $, $x_2 \in \mathbb{R}^m $. Suppose that
$\overline{\Omega}$ is invariant under the action of $G_{k,m}$
$(\tau (\overline{\Omega})=\overline{\Omega}, \forall \tau \in
G_{k,m})$. Then $\overline{\Omega}$ is a compact $n-$dimensional
manifold with boundary, invariant under the action of  the
subgroup $G_{k,m}$ of the isometry group $O(n)$.

\section{Notations and preliminary results}

$\,\,\,\,\,\,\,\,$ Let $(M,g)$ be a compact $n-$dimensional,
$n\geq 3$, Riemannian manifold with boundary $G-$invariant under
the action of a subgroup $G$ of the isometry group $I(M,g)$. We
assume that $(M,g)$ is a smooth bounded open subset of a slightly
larger Riemannian manifold $( {\widetilde M,g} ) $ (i.e. see
\cite{Li-Zhu}), invariant under the action of a subgroup $G$ of
the isometry group of $( {\widetilde M,g} ) $.\\
Consider the spaces of all $G-$invariant functions under the
action  of the group $G$
$$
C_G^\infty (M)=\left\{u\in C^\infty(M):u\circ \tau=u\, ,\, \forall
\: \tau \in G\right\}
$$
$$
C_{0,G}^\infty (M)=\left\{u\in C_0^\infty(M):u\circ \tau=u\, ,\,
\forall \: \tau \in G\right\}
$$
Denote $H_1^p(M)$ the completion of $C^\infty(M)$ with respect to
the norm
$$
\left\| u \right\|_{H_1^p \left( M \right)}  = \left( {\left\|
{\nabla u} \right\|_{L^p \left( M \right)}^p  + \left\| u
\right\|_{L^p \left( M \right)}^p } \right)^{1/p} ,
$$
and $H_{1, G}^p(M)$ the space of all $G-$invariant functions of
$H_1^p(M)$.\\
For completeness we cite some background material and results from \cite{Cot-Lab3}.\\
Given $(\widetilde M,g)$ a Riemannian manifold (complete or not,
but connected), we denote by $I(\widetilde M, g)$ its group of
isometries.

Let $ P \in M$ and $O_{ P}=\{\tau ( P), \tau \in G \} $ be its
orbit of dimension $k$, $0 \leqslant k < n$. According to
(\cite{Heb} \S \,9, \cite{Fag1}) the map $\Phi :G \to O _{ P} $,
defined by $\Phi \left( \tau \right) = \tau \left( { P} \right)$,
is of rank $k$ and there exists a submanifold $H$ of $G$ of
dimension $k$ with $Id \in H$, such that $\Phi $ restricted to $H$
is a diffeomorphism from $H$ onto its image denoted
$\mathcal{V}_{P } $.

Let $ N$ be a submanifold of $ M$ of dimension $(n-k)$, such that
$T_P \Phi \left( H \right) \oplus T_P  N = T_P  M$. Using the
exponential map at $ P$, we build a $\left( {n - k} \right)-$
dimensional submanifold $\mathcal{W}_{P} $ of $ N$, orthogonal to
$O_{ P} $ at ${P}$ and such that for any $ Q \in \mathcal{W}_{P}
$, the minimizing geodesics of $ \left( { M,g} \right)$ joining $
P$ and $ Q$ are all contained in $\mathcal{W}_{ P} $.

Let $\Psi:H \times \mathcal{W}_{ P} \to M$, be the map defined by
$\Psi \left( {\tau , Q} \right) = \tau \left( { Q} \right)$.
According to the local inverse theorem, there exists a
neighborhood $\mathcal{V}_{\left( {Id, P} \right)} \subset H
\times \mathcal{W}_{ P} $ of $\left( {Id, P} \right)$ and a
neighborhood $\mathcal{M}_{ P}  \subset M$ such that $\Psi ^{ - 1}
= \left({\Psi _1 \times \Psi _2 } \right)$, from $\mathcal{M}_{
P}$ onto $\mathcal{V}_{\left( {Id, P} \right)} $ is a
diffeomorphism.

Up to restricting $\mathcal{V}_{ P}$, we choose a normal chart
$\left( {\mathcal{V}_{ P} ,\varphi _1 } \right)$ around $ P$ for
the metric $\widetilde{g}$ induced on $O_{ P} $, with
$\varphi_1\left({\mathcal{V}_{ P} } \right) = U \subset
\mathbb{R}^k $. In the same way, we choose a geodesic normal chart
$\left( {\mathcal{W}_{ P} ,\varphi _2 } \right)$ around $ P$ for
the metric $\tilde{ \tilde {g}}$ induced on $\mathcal{W}_{ P} $,
with $\varphi _2 \left( {\mathcal{W}_{ P} } \right) = W \subset
\mathbb{R}^{n - k} $.\\
We denote by $\xi _1  = \varphi _1 \circ \Phi \circ \Psi _1 $,
$\xi _2  = \varphi _2 \circ  \Psi _2 $, $\xi = \left( {\xi _1 ,\xi
_2 } \right)$ and $\Omega = \mathcal{M}_{ P} $.

From the above and the Lemmas 1 and 2 in \cite{Heb-Vau} the
following lemma holds:
\begin{lemma}\label{L3.3}
Let $ \left( {M,g} \right)$ be a compact Riemannian $n-$manifold
with boundary, $G$ a compact subgroup of $I\left( { M,g} \right)$,
$ P \in  M$ with orbit of dimension $k$, $0 \leqslant k < n$. Then
there exists a chart $\left( {\Omega ,\xi } \right)$ around $P$
such that the following properties are valid:
\begin{enumerate}
    \item $\xi \left( \Omega  \right) = U \times W$, where $U
           \subset \mathbb{R}^k $ and $W \subset \mathbb{R}^{n - k} $.
    \item $U$, $W$ are bounded, and $W$ has smooth boundary.
    \item  $\left( {\Omega ,\xi } \right)$ is a normal chart of
          $ M$ around of $P$, $\left(
          {\mathcal{V}_{ P}, \varphi _1 } \right)$ is a normal
          chart around of $ P$ of submanifold $O_{P} $
          and $\left( {\mathcal{W}_{ P} ,\varphi _2 } \right)$ is
          a normal geodesic  chart around of $ P$ of submanifold
          $\mathcal{W}_{ P}$.
    \item For any $\varepsilon>0$, $\left( {\Omega ,\xi }
          \right)$ can be chosen such that:
          \[
          1 - \varepsilon \leq \sqrt {\det \left( {g_{ij} } \right)}  \leq 1 + \varepsilon
          \,\, \mathrm{on}\,\, \Omega , \,\,\mathrm{for}\; 1 \leq i,j \leq n
          \]
          \[
          1 - \varepsilon  \leq \sqrt {\det \left( {\tilde g_{ij} } \right)}
            \leq 1 + \varepsilon
          \,\,\mathrm{on} \,\,\mathcal{V}_{ P} , \,\,\mathrm{for
          }\; 1 \leq i,j \leq k.
          \]
    \item For any $u \in C_G^\infty(M)  $, $u \circ \xi ^{ - 1} $
    depends only on $W$ variables.
\end{enumerate}
\end{lemma}

We say that we choose a neighborhood of $O_P$ when we choose
$\delta>0$ and we consider
$$
 O_{\!P,\; \delta }  = \{ {Q \in
\widetilde M : d(Q,O_P ) < \delta }\}.
$$
Such a neighborhood of $O_P$ is called a tubular neighborhood.

Let $P \in M$ and $O_P $ be its orbit of dimension $k$. Since the
manifold $M$ is included in $\widetilde M$, we can choose a normal
chart $( {\Omega_P ,\xi_P } )$ around $P$ such that Lemma
\ref{L3.3} holds for some $\varepsilon_0>0$. For any $Q = \tau(P)
\in O_P $, where $\tau  \in G$, we build a chart around $Q$,
denoted by $( {\tau ( \Omega_P )\!, \xi_P \circ \tau ^{ - 1} } )$
and ``isometric'' to $( {\Omega_P ,\xi_P } )$. $ O_P $ is then
covered by such charts. We denote by $( {\Omega _{P,m} } )_ {m =
1,...,M} $ a finite extract covering. Then we can choose
$\delta>0$ small enough, depending on $P$ and $\varepsilon_0$ such
that the tubular neighborhood $ O_{\!P,\; \delta }$, (where $d( {
\cdot ,O_P }
)$ is the distance to the orbit) has the following properties:\\
$(i)\;$ $\overline { O_{\!P,\; \delta }} $ is a
submanifold of $\widetilde M$ with boundary,\\
$(ii)\;$ $d^2 ( { \cdot ,O_P } )$,  is a $C^\infty $
function on $ O_{\!P,\; \delta } $ and \\
$(iii)\;$ $ O_{\!P,\; \delta } $ is covered by $(
{\Omega _m } )_{m = 1,...,M} $.\\
Clearly, $M$ is covered by $\cup_{P \in M} O_{\!P,\; \delta } $.
We denote by $( {O_{j, \,\delta } } )_{j = 1,...,J} $ a finite
extract covering of $M$, where all $O_{j, \,\delta } $'s are
covered by $( {\Omega _{jm} } )_{m = 1,...,M_j } $. Then  we will
have
$$
M\subset \bigcup\nolimits_{j = 1}^J {\bigcup\nolimits_{m = 1}^{M_j
} {\Omega _{jm} } }=\bigcup\nolimits_{m = 1}^{\sum\nolimits_{j =
1}^J {M_j } }{\Omega _i }
$$
So we obtain a finite covering of $M$ consisting of $\Omega_i$'s,
$i = 1,...,\sum\nolimits_{j = 1}^J {M_j } $. We choose  such a
covering in the following way:

$(i)\;$ If $P$ lies in the interior of $M$, then there exist $j,\;
1\leq j\leq J$ and $m,\; 1\leq m\leq M_j$ such that the tubular
neighborhood $O _{j,\, \delta} $ and $\Omega _{jm} $, with
$P\in\Omega _{jm} $,  lie entirely  in $M$'s interior, (that is,
if $P \in M\backslash \partial M$, then $O_{j, \delta}\subset
M\backslash
\partial M $ and $\Omega_{jm} \subset M\backslash
\partial M$).

$(ii)\;$ If $P$ lies on the boundary  $\partial M$ of $M$, then a
$j,\; 1\leq j\leq J$ exists, such that the tubular neighborhood $O
_{j,\, \delta} $ intersects the boundary  $\partial M$ and an
$m,\; 1\leq m\leq M_j$ exists, such that $\Omega _{jm} $, with
$P\in\Omega _{jm} $, cuts a part of  the boundary  $\partial M$.
Then the $\Omega _{jm} $ covers a patch of the boundary of $M$,
and the whole of the  boundary is covered by charts around $P \in
\partial M$.\\

We denote $N$ the projection of the image of $M$, through the
charts $( {\Omega _{jm} ,\xi _{jm} } )$, ${j = 1,...,J } $, ${m =
1,...,M_j } $,  on $\mathbb{R}^{n - k} $. Then $( {N,\bar g} )$ is
a $(n-k)-$ dimensional compact submanifold of $\mathbb{R}^{n - k}
$ with boundary and $N$ is covered by $( {W_i } )$, $i =
1,...,\sum\nolimits_{j = 1}^J {M_j } $, where $W_i $ is the
component of $\xi _i ( {\Omega _i } )$ on $\mathbb{R}^{n - k} $
for all $i = 1,...,\sum\nolimits_{j = 1}^J
{M_j }$.\\
Let $p $  be the projection of $\xi _i ( {P } ),P \in M$ on
$\mathbb{R}^{n - k} $. Thus one of the following holds:

$(i)\;$ If $p  \in N\backslash \partial N$, then $W_i \subset
N\backslash
\partial N$ and $W_i$ is a normal geodesic neighborhood with
normal geodesic coordinates $( {y_1 ,...,y_{n - k} } )$.

$(ii)\;$ If $p  \in \partial N$, then $W_i$ is a Fermi
neighborhood with Fermi coordinates $( {y_1 ,...,y_{n - k - 1} ,t} )$.\\
In these neighborhoods we have
 \[1 - \varepsilon _0  \leqslant \sqrt {\det \left( {\bar
g_{ij} } \right)}  \leqslant 1 + \varepsilon _0\,\,
\mathrm{on}\,\,N,\,\,\mathrm{for}\,\,1\leq i,j\leq n-k,\] where
$\varepsilon_0$ can be as small as we want, depending on the
chosen covering.\\

For convenience in the following we set
$$
O_j  =  O_{j,\, \delta} = \{ {Q \in \widetilde M : d(Q,O_{P  _j }
) < \delta } \}
$$
We still need the following  lemma:
\begin{lemma}\label{L3.4}
$\emph{\textbf{(a)}\,\,\,}$ For any $\upsilon  \in H_{1,G}^p
\left( {O_j \cap M} \right),\upsilon  \geqslant 0$ the following
properties are valid:
\begin{enumerate}
    \item  $$ \left( {1 - c\varepsilon _0 } \right)V_j
\int_N { \upsilon _2 ^p } d\upsilon _{\bar g}   \leqslant \int_M {
\upsilon } ^p dV_g \leqslant \left( {1 + c\varepsilon _0 }
\right)V_j \int_N { \upsilon _2 ^p} d\upsilon _{\bar g} $$
    \item  $$ \left( {1 - c\varepsilon _0 } \right)V_j
\int_N {\left| {\nabla _{\bar g} \upsilon _2 } \right|^p d\upsilon
_{\bar g} }  \leqslant  \int_M {\left| {\nabla _g \upsilon }
\right|^p dV_g }\leqslant  \left( {1 + c\varepsilon _0 }
\right)V_j \int_N {\left| {\nabla _{\bar g} \upsilon _2 }
\right|^p d\upsilon _{\bar g} }$$
\end{enumerate}
$\emph{\textbf{(b)}\,\,\,}$ For any $\upsilon  \in H_{1,G}^p
\left( {O_j  \cap \partial M} \right),\upsilon  \geqslant 0$ the
following property is valid:
$$ \left( {1 - c\varepsilon _0 } \right)V_j
\int_{\partial N} {\upsilon _2 ds_{\bar g} }  \leqslant
\int_{\partial M} {\upsilon dS_g }  \leqslant \left( {1 +
c\varepsilon _0 } \right)V_j \int_{\partial N} {\upsilon _2
ds_{\bar g} }$$ where $V_j  = Vol\left( {O_j } \right)$,
$\upsilon_2=\upsilon\circ\xi^{-1}$ and $c$ is a positive constant.
\end{lemma}

Moreover, we need the following propositions:
\begin{proposition}\label{P4.1}
For any $\varepsilon  > 0$ and for all $ u \in C_G^\infty  ( M ) $
the following inequality holds
\begin{eqnarray}\label{E4.1}
\left( {\int_{\partial M} {\left( {\eta _j u} \right)^2 dS_g } }
\right)^{\frac{{n - k}} {{n - k - 1}}} & \leqslant& \frac{{\tilde
A_0 \left( {n - k} \right) + \varepsilon }} {{V_j^{\frac{1} {{n -
k - 1}}} }}\int_M {\left| {\nabla \left( {\eta _j u} \right)}
\right|^2 dV_g }\nonumber \\&&\times\left( {\int_{\partial M}
{\eta _j\left| u \right|dS_g } } \right)^{\frac{2} {{n - k - 1}}}
\end{eqnarray}
where $(\eta_j)$ is a partition of unity associating to $(O_j)$.
\end{proposition}
{\bf{Proof}}. By Lemma \ref{L3.4} (b) with $\upsilon  =  \eta _j u
$ and $p=2$ we obtain
\begin{equation}\label{E4.2}
\int_{\partial M} {\left( {\eta _j u} \right)^2 dS_g } \,\,\,\,
\leqslant \left( {1 + c\varepsilon _0 } \right)V_j \int_{\partial
N} {\left( {\eta _j u} \right)_2^2 ds_{\bar g} }
\end{equation}
Let $\varepsilon _0 > 0$. Then there exist $\delta>0$ such that
for any $Q = \xi( P ) \in \partial N,\, P \in \partial M$ and for
all $\phi  \in C_0^\infty  ( {B_Q ( \delta  )} )$,\, ($B_Q( \delta
)$ is the $(n-k)-$dimensional ball of radius $\delta$ centered on
$Q$), according to Theorem 2 of \cite{Hum2} the following
inequality holds
\begin{equation}\label{E4.3}
\left( {\int_{\partial N} {\phi ^2 ds_{\bar g} } }
\right)^{\frac{{n - k}} {{n - k - 1}}}  \leqslant \left( {\tilde
A_0 \left( {n - k} \right) + \frac{\varepsilon_0 } {2}}
\right)\int_N {\left| {\nabla _{\bar g} \phi } \right|^2 }
d\upsilon_{\bar g} \left( {\int_{\partial N} {\left| \phi
\right|ds_{\bar g} } } \right)^{\frac{2} {{n - k - 1}}}
\end{equation}
where $\tilde A_0 ( {n - k} )$ is the best constant of the trace
Nash inequality
\begin{eqnarray*}
\left( {\int_{\partial N} {\phi ^2 ds_{\bar g} } }
\right)^{\frac{{n - k}} {{n - k - 1}}}  \leqslant \left( {\tilde
A\int_N {\left| {\nabla \phi } \right|^2 } d\upsilon_{\bar g}  +
\tilde B\int_{\partial N} {\phi ^2 ds_{\bar g} } } \right)\left(
{\int_{\partial N} {\left| \phi \right|ds_{\bar g} } }
\right)^{\frac{2} {{n - k - 1}}}
\end{eqnarray*}
By (\ref{E4.2}) and (\ref{E4.3}) we have
\begin{eqnarray}\label{E4.5}
\left( {\int_{\partial M} {\left( {\eta _j u} \right)^2 dS_g }
\,\,} \right)\,^{\frac{{n - k}} {{n - k - 1}}} \, &\leqslant&
\left[ {\left( {1 + c\varepsilon _0 } \right)V_j }
\right]^{\frac{{n - k}} {{n - k - 1}}} \left( {\tilde A_0 \left(
{n - k} \right) + \frac{\varepsilon_0 }
{2}} \right) \nonumber \\
&& \times \int_N {\left| {\nabla _{\bar g} \left( {\eta _j u}
\right)_2 } \right|^2 } d\upsilon_{\bar g} \left( {\int_{\partial
N} {\left| {\left( {\eta _j u} \right)_2 } \right|ds_{\bar g} } }
\right)^{\frac{2} {{n - k - 1}}} \nonumber \\
\end{eqnarray}
From (\ref{E4.5}) and Lemma \ref{L3.4} arises
\begin{eqnarray*}
\left( {\int_{\partial M} {\left( {\eta _j u} \right)^2 dS_g }
\,\,} \right)\,^{\frac{{n - k}} {{n - k - 1}}} \, &\leqslant &
\frac{{\left[ {\left( {1 + c\varepsilon _0 } \right)V_j }
\right]^{\frac{{n - k}} {{n - k - 1}}} }} {{\left[ {\left( {1 -
c\varepsilon _0 } \right)V_j } \right]^{\frac{2} {{n - k - 1}}}
\left( {1 - c\varepsilon _0 } \right)V_j }}\\ && \times\left(
{\tilde A_0 \left( {n - k} \right) + \frac{\varepsilon_0 } {2}}
\right) \int_M {\left| {\nabla _{ g} \left( {\eta _j u} \right)}
\right|^2 } dV_g \\ && \times\left( {\int_{\partial M} {\left|
{\left( {\eta _j u} \right)} \right|dS_g } } \right)^{\frac{2} {{n
- k - 1}}}
\end{eqnarray*}
or
\begin{eqnarray}\label{E4.6}
\left( {\int_{\partial M} {\left( {\eta _j u} \right)^2 dS_g }
\,\,} \right)\,^{\frac{{n - k}} {{n - k - 1}}} \, &\leqslant&
\frac{{\left( { {1 + c\varepsilon _0 }  } \right)^{\frac{{n - k}}
{{n - k - 1}}} }} {{ {\left( {1 - c\varepsilon _0 } \right) }
^{\frac{n-k+1} {{n - k - 1}}}  }}\left( {\tilde A_0 \left( {n - k}
\right) + \frac{{\varepsilon _0 }} {2}} \right) \frac{1}
{{V_j^{\frac{1} {{n - k - 1}}} }} \nonumber
\\ && \times\int_M {\left| {\nabla _{ g}
\left( {\eta _j u} \right)} \right|^2 } dV_g \left(
{\int_{\partial M} {\left| {\left( {\eta _j u} \right)}
\right|dS_g } } \right)^{\frac{2} {{n - k - 1}}}\nonumber\\
\end{eqnarray}

Given $\varepsilon>0$ we can choose  $\varepsilon_0>0$ small
enough such that
$$\frac{{\left( { {1 + c\varepsilon _0 }  } \right)^{\frac{{n - k}}
{{n - k - 1}}} }} {{ {\left( {1 - c\varepsilon _0 } \right) }
^{\frac{n-k+1} {{n - k - 1}}}  }}\left( {\tilde A_0 \left( {n - k}
\right) + \frac{{\varepsilon _0 }} {2}} \right) \leqslant \tilde
A_0 ( {n - k} ) + \varepsilon
$$ and then by (\ref{E4.6}) we obtain (\ref{E4.1}).
\begin{proposition}\label{P4.2}
For any $\varepsilon  > 0$ and for all $ u \in C_G^\infty ( M )$
the following inequality holds
\begin{eqnarray}\label{E4.7}
  \left( {\int_{\partial M} {u^2 dS_g } } \right)^{\frac{{n - k}}
{{n - k - 1}}} & \leqslant & {\frac{{\tilde A_0\left( {n - k}
\right) + \varepsilon }} {{V^{\frac{1} {{n - k - 1}}} }}} \left(
{\int_M {\left| {\nabla u} \right|^2 dV_g } + C \int_M {u^2 dV_g }
} \right)\nonumber\\&&\times \left( {\int_{\partial M} {\left| u
\right|} dS_g } \right)^{\frac{2} {{n - k - 1}}}
\end{eqnarray}
\end{proposition}
{\bf{Proof}}. We set $\alpha _j = \frac{{\eta _j^2 }}
{{\sum\nolimits_{m = 1}^{M_j} {\eta _m^2 } }}, j = 1,2,..,J$ and
so $ \{ {\alpha _j } \} $ is a partition of unity for $ M $
subordinated in the covering $( {O_j })_{j = 1,2,...,J}$ and
functions $\sqrt {\alpha _j } $ are smooth, $G-$invariants and
there exist a constant $H$ such that for any $j = 1,...,J$ holds
\begin{equation}\label{E4.8}
| {\nabla \sqrt {\alpha _j } } | \leqslant H
\end{equation}
Let $ u \in C_G^\infty ( M )$. Then we have
\begin{equation}\label{E4.9}
\int_{\partial M} {u^2 dS_g  = } \int_{\partial M} {\left(
{\sum\nolimits_{j = 1}^J {\alpha _j } u^2 } \right)dS_g  = }
\sum\nolimits_{j = 1}^J {\int_{\partial M} {\left( {\sqrt {\alpha
_j } u} \right)^2 dS_g } }
\end{equation}
By (\ref{E4.9}) and Proposition \ref{P4.1} arises
\begin{eqnarray*}
  \int_{\partial M} {u^2 dS_g }& \leqslant &  \sum\limits_{j =
  1}^J
  {\left( {\frac{{\tilde A_0\left( {n - k} \right) + \varepsilon }}
{{V_j^{\frac{1} {{n - k - 1}}} }}} \right)^{\frac{{n - k - 1}} {{n
- k}}} \left( {\int_M {\left| {\nabla \left( {\sqrt {\alpha _j }
u} \right)} \right|^2 dV_g } } \right)^{\frac{{n - k - 1}} {{n -
k}}} } \\&&\times \left( {\int_{\partial M} {\sqrt {\alpha _j }
\left| u \right|} dS_g } \right)^{\frac{2} {{n - k}}}
\end{eqnarray*}
and since $\min V_j  = V$ we obtain
\begin{eqnarray}\label{E4.10}
\int_{\partial M} {u^2 dS_g  }&\leqslant & \left( {\frac{{\tilde
A_0\left( {n - k} \right) + \varepsilon }} {{V^{\frac{1} {{n - k -
1}}} }}} \right)^{\frac{{n - k - 1}} {{n - k}}} \sum\limits_{j =
1}^J {\left( {\int_M {\left| {\nabla \left( {\sqrt {\alpha _j } u}
\right)} \right|^2 dV_g } } \right)^{\frac{{n - k - 1}} {{n - k}}}
}\nonumber\\&&\times \left( {\int_{\partial M} {\sqrt {\alpha _j }
\left| u \right|} dS_g } \right)^{\frac{2} {{n - k}}}
\end{eqnarray}
By H\"older's inequality we have
\begin{eqnarray*}
  \left( {\int_{\partial M} {\sqrt {\alpha _j } \left| u \right|} dS_g } \right)^{\frac{2}
{{n - k}}} & = &\left[ {\int_{\partial M} {\left( {\sqrt {\alpha
_j \left| u \right|} \sqrt {\left| u \right|} } \right)} dS_g }
\right]^{\frac{2} {{n - k}}} \\ & \leqslant & \left(
{\int_{\partial M} {\alpha _j \left| u \right|} dS_g }
\right)^{\frac{1} {{n - k}}} \left( {\int_{\partial M} {\left| u
\right|} dS_g } \right)^{\frac{1} {{n - k}}}
\end{eqnarray*}
and by (\ref{E4.10}) we obtain
\begin{eqnarray}\label{E4.11}
\int_{\partial M} {u^2 dS_g  }&\leqslant &  \left( {\frac{{\tilde
A_0\left( {n - k} \right) + \varepsilon }} {{V^{\frac{1} {{n - k -
1}}} }}} \right)^{\frac{{n - k - 1}} {{n - k}}} \left(
{\int_{\partial M} {\left| u \right|} dS_g } \right)^{\frac{1} {{n
- k}}}\nonumber \\
 && \times   \sum\limits_{j = 1}^J {\left( {\int_M {\left| {\nabla
\left( {\sqrt {\alpha _j } u} \right)} \right|^2 dV_g } }
\right)^{\frac{{n - k - 1}} {{n - k}}} } \left( {\int_{\partial M}
{\alpha _j \left| u \right|} dS_g } \right)^{\frac{1} {{n -
k}}}\nonumber\\
\end{eqnarray}
Moreover the following H\"older's inequality
\begin{equation}\label{E4.12}
\sum\nolimits_{j = 1}^J {a_j b_j } \leqslant \left(
{\sum\nolimits_{j = 1}^J {a_j^p } } \right)^{1/p} \left(
{\sum\nolimits_{j = 1}^J {b_j^q } } \right)^{1/q}
\end{equation}
 holds for any $a_j ,b_j $ nonnegative and for all $p \geqslant 1,q
\geqslant 1$ with $( {1/p} ) + ( {1/q} ) =
1$.\\
Setting  in (\ref{E4.12}) $$ a_j  = \left( {\int_M {| {\nabla (
{\sqrt {\alpha _j } u} )} |^2 dV_g } }\right)^{\frac{{n - k - 1}}
{{n - k}}}, \quad b_j  = \left( {\int_{\partial M} {\alpha _j | u
|} dS_g} \right)^{\frac{1} {{n - k}}}$$$$p = \frac{{n - k}} {{n -
k - 1}}, \quad q = n - k$$ we obtain
\begin{eqnarray}\label{E4.13}
&& \sum\limits_{j = 1}^J {\left( {\int_M {\left| {\nabla \left(
{\sqrt {\alpha _j } u} \right)} \right|^2 dV_g } }
\right)^{\frac{{n - k - 1}} {{n - k}}} } \left( {\int_{\partial M}
{\alpha _j \left| u \right|} dS_g } \right)^{\frac{1} {{n - k}}}\nonumber \hfill \\
&\leqslant& \left( {\sum\limits_{j = 1}^J {\int_M {\left| {\nabla
\left( {\sqrt {\alpha _j } u} \right)} \right|^2 dV_g } } }
\right)^{\frac{{n - k - 1}} {{n - k}}} \left( {\sum\limits_{j =
1}^J {\int_{\partial M} {\alpha _j \left| u \right|} dS_g } }
\right)^{\frac{1} {{n - k}}} \nonumber\hfill \\ & = &\left(
{\sum\limits_{j = 1}^J {\int_M {\left| {\nabla \left( {\sqrt
{\alpha _j } u} \right)} \right|^2 dV_g } } } \right)^{\frac{{n -
k - 1}} {{n - k}}} \left( {\int_{\partial M} {\left(
{\sum\limits_{j = 1}^J {\alpha _j } } \right)\left| u \right|dS_g
} } \right)^{\frac{1} {{n - k}}}\nonumber \hfill \\ &  =&
  \left( {\int_{\partial M} {\left| u \right|dS_g } } \right)^{\frac{1}
{{n - k}}} \left( {\sum\limits_{j = 1}^J {\int_M {\left| {\nabla
\left( {\sqrt {\alpha _j } u} \right)} \right|^2 dV_g } } }
\right)^{\frac{{n - k - 1}} {{n - k}}}
\end{eqnarray}
By (\ref{E4.11}) and (\ref{E4.13}) arises
\begin{eqnarray}\label{E4.14}
\int_{\partial M} {u^2 dS_g } &\leqslant & \left( {\frac{{\tilde
A_0\left( {n - k} \right) + \varepsilon }} {{V^{\frac{1} {{n - k -
1}}} }}} \right)^{\frac{{n - k - 1}} {{n - k}}} \left(
{\int_{\partial M} {\left| u \right|} dS_g } \right)^{\frac{2} {{n
- k}}}  \nonumber\\ &&\times \left( {\sum\limits_{j = 1}^J {\int_M
{\left| {\nabla \left( {\sqrt {\alpha _j } u} \right)} \right|^2
dV_g } } } \right)^{\frac{{n - k - 1}} {{n - k}}}
\end{eqnarray}
Further more since
\[
\left| {\nabla \left( {\sqrt {\alpha _j } u} \right)} \right|^2  =
\alpha _j \left| {\nabla u} \right|^2  + u^2 \left| {\nabla \left(
{\sqrt {\alpha _j } } \right)} \right|^2 + 2\left\langle {\nabla
u,\nabla \left( {\sqrt {\alpha _j } } \right)} \right\rangle
u\sqrt {\alpha _j }
\]
and since (\ref{E4.8}) holds, after some computations  we obtain
\begin{eqnarray*}
\sum\limits_{j = 1}^J {\int_M {\left| {\nabla  \left( {\sqrt
{\alpha _j } u} \right)} \right|^2 dV_g }}
 & \leqslant & \int_M {\left|
  {\nabla u} \right|^2 dV_g }  + HJ\int_M
  {u^2 dV_g }  \hfill \\
 && + 2\int_M
  {\sum\limits_{j = 1}^J {\left\langle {\nabla
   \left( {\sqrt {\alpha _j } } \right),\nabla u}
   \right\rangle u\sqrt {\alpha _j } dV_g } }
\end{eqnarray*}
But
\begin{eqnarray*}
0&=&\nabla ( {\sum\nolimits_{j = 1}^J {\alpha _j } })
 =\sum\nolimits_{j = 1}^J {\left( {\nabla \alpha _j } \right)}\\
 &=&\sum\nolimits_{j = 1}^J {2\sqrt {\alpha _j } \nabla \left( {\sqrt
  {\alpha _j } } \right)}=2\sum\nolimits_{j = 1}^J {\sqrt {\alpha _j}
  \nabla \left( {\sqrt {\alpha _j } } \right) }
\end{eqnarray*}
Thus the following inequality holds
\begin{equation}\label{E4.15}
\sum\limits_{j = 1}^J {\int_M {\left| {\nabla \left( {\sqrt
{\alpha _j } u} \right)} \right|^2 dV_g } }  \leqslant \int_M
{\left| {\nabla u} \right|^2 dV_g }  + C\int_M {u^2 dV_g }
\end{equation}
Finally by (\ref{E4.14}) and (\ref{E4.15}) we have
\begin{eqnarray*}
\int_{\partial M} {u^2 dS_g }& \leqslant & \left( {\frac{{\tilde
A_0\left( {n - k} \right) + \varepsilon }} {{V^{\frac{1} {{n - k -
1}}} }}} \right)^{\frac{{n - k - 1}} {{n - k}}} \left(
{\int_{\partial M} {\left| u \right|} dS_g } \right)^{\!\frac{2}
{{n - k}}} \\ &&\times\left( {\int_M {\left| {\nabla u} \right|^2
dV_g }\! +\! C\int_M {u^2 dV_g } } \right)^{\frac{{n - k - 1}}
{{n - k}}}  \hfill \\
\end{eqnarray*}
and the proposition is proved.

\section{Proofs}

\noindent{\bf{Proof of Theorem \ref{T2.1}}}. The proof is based on
\cite{Dru-Heb-Vau}. Let as sketch the proof. Following the same
steps as in Proposition \ref{P4.1} as well as by Lemma \ref{L3.4}
and Theorem $1.1$ in \cite{Dru-Heb-Vau} we obtain  that for any
$\varepsilon_0  > 0$ and for all $ u \in C_G^\infty  ( M ) $ the
following inequality holds
\begin{eqnarray}\label{E4.21}
\left( {\int_{ M} {\left( {\eta _j u} \right)^2 dV_g } }
\right)^{\frac{{n - k +2}} {{n - k }}} & \leqslant & \frac{{ A_0
\left( {n - k} \right) + \frac{\varepsilon_0}{2} }}
{{V_j^{\frac{2} {{n - k }}} }}\int_M {\left| {\nabla \left( {\eta
_j u} \right)} \right|^2 dV_g }\nonumber\\&&\times \left( {\int_{
M} {\eta _j\left| u \right|dV_g } } \right)^{\frac{4} {{n - k }}},
\end{eqnarray}
where $(\eta_j)$ is a partition of unity associating to $(O_j)$.\\
Let $\{\alpha _j\}_{ j = 1,2,..,J}$  a partition of unity for $ M
$ as in Proposition \ref{P4.2}. Since for any $ u \in C^\infty ( M
) $ the following holds
$$ \|u\|^2_2=\|u^2\|_1=\|\Sigma_{j = 1}^J \alpha_j u^2\|_1\leq
\Sigma_{j = 1}^J\|\alpha _j u^2\|_1=\Sigma_{j = 1}^J\|\sqrt{\alpha
_j} u\|^2_2$$ where $\|.\|$ stands for the $L^p-$norm.\\ By
H\"older's inequality $\|\sqrt{\alpha _j}u\|_1\leq \|\alpha_j
u\|_1^{1/2}\|u\|_1^{1/2} $ from (\ref{E4.21}) we obtain
\begin{eqnarray}\label{E4.22}
\int_M {u^2 dV_g }   &\leqslant& \left( \frac{ A_0 \left( {n - k}
\right) +  \frac{\varepsilon_0}{2}} {V^{\frac{2} {{n - k }}}
}\right)^\frac{n-k}{n-k+2}\left( {\int_{ M} {\left| u \right|dV_g
} }\right)^{\frac{2} {{n - k+2 }}}\nonumber\\
&&\times\sum_{j=1}^J\left(\int_M {\left| {\nabla \left( {\alpha _j
u} \right)} \right|^2 dV_g }\right)^\frac{n-k}{n-k+2} \left(
{\int_{ M} {\alpha _j\left|  u \right|dV_g } } \right)^{\frac{2}
{{n - k +2}}}\nonumber\\
\end{eqnarray}
By (\ref{E4.22}) because of (\ref{E4.8}) and (\ref{E4.12}) with
$$
 a_j  = \left( {\int_M {| {\nabla ( {\sqrt {\alpha _j } u} )}
|^2 dV_g } }\right)^{\frac{{n - k}} {{n - k+2}}}, \quad b_j  =
\left( {\int_{\partial M} {\alpha _j | u |} dV_g}
\right)^{\frac{2} {{n - k+2}}}
$$
$$
p = \frac{n - k+2} {n - k },\quad q = \frac{n - k+2}{2}
$$
arises
\begin{eqnarray}\label{E4.22}
\left( {\int_M {u^2 dV_g } } \right)^{\frac{{n - k + 2}} {{n -
k}}} & \leqslant & \left( \frac{ A_0 \left( {n - k} \right) +
\frac{\varepsilon_0}{2}} {V^{\frac{2} {{n - k }}}
}\right)^\frac{n-k}{n-k+2} \nonumber\\&&\times\left({\int_M
{\left| {\nabla u} \right|_g^2 } dV_g  + H^2J \int_M {u^2 dV_g } }
\right) \nonumber\\&&\times \left( {\int_M {\left| u \right|dV_g }
} \right)^{\frac{4} {{n - k}}}
\end{eqnarray}
Given $\varepsilon>0$ we can choose $\varepsilon_0>0$ such that
$$
 \frac{ A_0 \left( {n - k} \right) +
\frac{\varepsilon_0}{2}} {V^{\frac{2} {{n - k }}} }\leq \frac{ A_0
\left( {n - k} \right)} {V^{\frac{2} {{n - k }}} } + \varepsilon
$$
and the theorem is proved.\\

\noindent{\bf{Proof of Theorem \ref{T2.2}}}. According to
\cite{Fag2} (Lemma 4) there exists an orbit of minimum dimension
$k$ and of minimum volume. Let $\mathcal{O}$ be the orbit of
dimension $k$ and of minimum volume, that is $Vol( \mathcal{O} ) =
\min Vol( {O_j } )_{j = 1,...,J} = V$. Let also the set $
\mathcal{O}_\delta = \{ {Q \in \widetilde M : d(Q,\mathcal{O}) <
\delta } \}$, where $ d( { \cdot ,\mathcal{O}} )$ is the distance
to the orbit. For $ u \in C_G^\infty ( {\mathcal{O}_\delta   \cap
M} ) $ by Proposition \ref{P4.2} because of Lemma \ref{L3.4} we
have sequentially

$$
\left( {\left( {1 - c\varepsilon _0 } \right)V\int_{\partial N}
{u_2^2 ds_{\bar g}}} \right)^{\frac{{n - k}} {{n - k - 1}}}
\leqslant\left( {\int_{\partial M} {u^2 dS_g } } \right)^{\frac{{n
- k}} {{n - k - 1}}}$$
\begin{eqnarray*} &\leqslant & \frac{{\tilde A_0\left( {n - k} \right) +
\frac{\varepsilon_0 } {2}}} {{V^{\frac{1} {{n - k - 1}}} }}\left(
{\int_M {\left| {\nabla u} \right|^2 dV_g } + C\int_M {u^2 dV_g }
} \right) \left( {\int_{\partial M} {\left| u \right|} dS_g }
\right)^{\frac{2} {{n - k - 1}}}\\
 & \leqslant &
\frac{{\tilde A_0\left( {n - k} \right) + \frac{\varepsilon_0 }
{2}}} {{V^{\frac{1} {{n - k - 1}}} }} \left( {\left( {1 +
c\varepsilon _0 } \right)V}{\int_N {\left| {\nabla u_2 } \right|^2
d\upsilon_{\bar g} }  + C{\left( {1 + c\varepsilon _0 } \right)V}
\int_N {u_2^2 d\upsilon_{\bar g}
} } \right)\hfill \\
&& \times\left( {\left( {1 + c\varepsilon _0 }
\right)V}{\int_{\partial N} {\left| {u_2 } \right|} ds_{\bar
g} } \right)^{\frac{2} {{n - k - 1}}}\\
 & \leqslant &
\frac{{\tilde A_0\left( {n - k} \right) + \frac{\varepsilon_0 }
{2}}} {{V^{\frac{1} {{n - k - 1}}} }}\left[ {\left( {1 +
c\varepsilon _0 } \right)V} \right]^{1 + {\frac{2} {{n - k - 1}}}}
\left( {\int_N {\left| {\nabla u_2 } \right|^2 d\upsilon_{\bar g}
}  + C  \int_N {u_2^2 d\upsilon_{\bar g} } }
\right)\\
&& \times\left( {\int_{\partial N} {\left| {u_2 } \right|}
ds_{\bar g} } \right)^{\frac{2} {{n - k - 1}}}
\end{eqnarray*}
or
\begin{eqnarray*}
  \left( {\int_{\partial N} {u_2^2 ds_{\bar g} } } \right)^{\frac{{n - k}}
{{n - k - 1}}} & \leqslant &\frac{{\tilde A_0\left( {n - k}
\right) + \frac{\varepsilon_0 } {2}}} {{V^{\frac{1} {{n - k - 1}}}
}}\frac{{\left( {1 + c\varepsilon _0 } \right)^{ \frac{n-k+1} {{n
- k - 1}}} }} {{\left( {1 - c\varepsilon _0 } \right)^{\frac{{n -
k}} {{n - k - 1}}} }}V^{1 + \frac{2} {{n - k - 1}} - \frac{{n -
k}}{{n - k - 1}}}  \hfill \\
&& \times \,\left( {\int_N {\left| {\nabla u_2 } \right|^2
d\upsilon_{\bar g} }  + C  \int_N {u_2^2 d\upsilon_{\bar g} } }
\right)\hfill \\
&& \times \left( {\int_{\partial N} {\left| {u_2 } \right|}
ds_{\bar g} } \right)^{\frac{2} {{n - k - 1}}},
\end{eqnarray*}
\begin{eqnarray*}
  \left( {\int_{\partial N} {u_2^2 ds_{\bar g} } } \right)^{\frac{{n - k}}
{{n - k - 1}}} & \leqslant & \frac{{\left( {1 + c\varepsilon _0 }
\right)^{\frac{n-k+1} {{n - k - 1}}} }} {{\left( {1 - c\varepsilon
_0 } \right)^{\frac{{n - k}} {{n - k - 1}}} }}\frac{{\tilde
A_0\left( {n - k} \right) + \frac{\varepsilon_0 } {2}}}
{{V^{\frac{1} {{n - k - 1}}} }}V^{\frac{1}
{{n - k - 1}}}  \hfill \\
&& \times \,\left( {\int_N {\left| {\nabla u_2 } \right|^2
d\upsilon_{\bar g} }  + C  \int_N {u_2^2 d\upsilon_{\bar g} } }
\right)\hfill \\
&& \times \left( {\int_{\partial N} {\left| {u_2 } \right|}
ds_{\bar g} } \right)^{\frac{2} {{n - k - 1}}} ,
\end{eqnarray*}
\begin{eqnarray*}
  \left( {\int_{\partial N} {u_2^2 ds_{\bar g} } } \right)^{\frac{{n - k}}
{{n - k - 1}}} & \leqslant &\frac{{\left( {1 + c\varepsilon _0 }
\right)^{\frac{n-k+1} {{n - k - 1}}} }} {{\left( {1 - c\varepsilon
_0 } \right)^{\frac{{n - k}} {{n - k - 1}}} }}\left( {\tilde
A_0\left( {n - k} \right) + \frac{\varepsilon_0 }
{2}} \right) \hfill \\
&& \times \,\left( {\int_N {\left| {\nabla u_2 } \right|^2
d\upsilon_{\bar g} }  + C  \int_N {u_2^2 d\upsilon_{\bar g} } }
\right)\hfill \\
&& \times \left( {\int_{\partial N} {\left| {u_2 } \right|}
ds_{\bar g} } \right)^{\frac{2} {{n - k - 1}}}
\end{eqnarray*}
hence, given $\varepsilon >0$, we choose $\varepsilon _0  > 0$
small enough to imply
\begin{eqnarray}\label{E4.16}
  \left( {\int_{\partial N} {u_2^2 ds_{\bar g} } } \right)^{\frac{{n - k}}
{{n - k - 1}}} & \leqslant & \left( {\tilde A_0\left( {n - k}
\right) + \varepsilon } \right)\,\left( {\int_N {\left| {\nabla
u_2 } \right|^2 d\upsilon_{\bar g} }
+ C  \int_N {u_2^2 d\upsilon_{\bar g} } } \right) \nonumber \\
&&\times \left( {\int_{\partial N} {\left| {u_2 } \right|}
ds_{\bar g} } \right)^{\frac{2} {{n - k - 1}}}
\end{eqnarray}
By (\ref{E4.16}) and Theorem 2 in \cite{Hum2} arises
\begin{eqnarray}\label{E4.17}
  \left( {\int_{\partial N} {u_2^2 ds_{\bar g} } } \right)^{\frac{{n - k}}
{{n - k - 1}}} & \leqslant &\,\left( {\left( {\tilde A_0\left( {n
- k} \right) + \varepsilon }
 \right)\int_N {\left| {\nabla u_2 } \right|^2 d\upsilon_{\bar g} }
  +\tilde B_\varepsilon  \int_{\partial N} {u_2^2 ds_{\bar g} } } \right) \nonumber \\
&& \times \left( {\int_{\partial N} {\left| {u_2 } \right|}
ds_{\bar g} } \right)^{\frac{2} {{n - k - 1}}}
\end{eqnarray}
Since the best constant of the  trace Nash inequality in $M$ has
the same value with the best constant of the trace Nash inequality
of the manifold $\mathcal{O}_\delta$, suppose that for any
$\alpha>0$ there exists $\theta > 0$ such that
\begin{equation}\label{E4.18}
\lambda _\alpha   = \mathop {\inf }\limits_A I_\alpha   \leqslant
\frac{{V^{\frac{1} {{n - k - 1}}} }} {{\tilde A_0( {n - k}) +
\varepsilon }} - \theta
\end{equation}
where
\[
I_\alpha   = \frac{{\left( {\int_M {\left| {\nabla u_\alpha }
\right|^2 dV_g }  + \alpha \int_{\partial M} {u_\alpha ^2 dS_g } }
\right)\left( {\int_{\partial M} {\left| {u_\alpha  } \right|}
dS_g} \right)^{\frac{2} {{n - k - 1}}} }} {{\left( {\int_{\partial
M} {u_\alpha ^2 dS_g } } \right)^{\frac{{n - k}} {{n - k - 1}}} }}
\]
and
\[
A = \{ {u \in H^2_{1,G} ( M ):u\left| {_{\partial M} \not \equiv
0} \right.} \}
\]
Thus there exists $u_\alpha   \in H^2_{1,G} ( {\mathcal{O}_\delta
} )$ such that
\begin{equation}\label{E4.19}
I_\alpha   \leqslant \frac{{V^{\frac{1} {{n - k - 1}}} }} {{\tilde
A_0( {n - k} ) + \varepsilon }} - \theta
\end{equation}
By (\ref{E4.19}) and Proposition \ref{P4.2},  we have
\begin{eqnarray*}
\frac{{\left[ {\left( {1 - c\varepsilon _0 } \right)V} \right]^{
\frac{n - k + 1} {{n - k - 1}}} \left( {\int_N {\left| {\nabla
_{\bar g} u_\alpha } \right|^2 d\upsilon _{\bar g} }  + \alpha
\int_{\partial N} {u_\alpha ^2 ds_{\bar g} } } \right)\left(
{\int_{\partial N} {\left| {u_\alpha  } \right|} ds_{\bar g} }
\right)^{\frac{2} {{n - k - 1}}} }} {{\left[ {\left( {1 +
c\varepsilon _0 } \right)V} \right]^{\frac{{n - k}} {{n - k - 1}}}
\left( {\int_{\partial N} {u_\alpha ^2 ds_{\bar g} } }
\right)^{\frac{{n - k}} {{n - k - 1}}} }} \leqslant\\
 \frac{{V^{\frac{1}
{{n - k - 1}}} }} {{\tilde A_0( {n - k} ) + \varepsilon }} -
\theta
\end{eqnarray*}
or
\begin{eqnarray*}
\frac{{\left( {\int_N {\left| {\nabla _{\bar g} u_\alpha }
\right|^2 d\upsilon _{\bar g} }  + \alpha \int_{\partial N}
{u_\alpha ^2 ds_{\bar g} } } \right)\left( {\int_{\partial N}
{\left| {u_\alpha } \right|} ds_{\bar g} } \right)^{\frac{2} {{n -
k - 1}}} }} {{\left( {\int_{\partial N} {u_\alpha ^2 ds_{\bar g} }
} \right)^{\frac{{n - k}} {{n - k - 1}}} }}
\end{eqnarray*}
\begin{eqnarray*}
& \leqslant&  \frac{{\left[ {\left( {1 + c\varepsilon _0 }
\right)V} \right]^{\frac{{n - k}} {{n - k - 1}}} }} {{\left[
{\left( {1 - c\varepsilon _0 } \right)V} \right]^{ \frac{n-k+1}
{{n - k - 1}}} }}\left[ {\left( {\frac{{\tilde A_0\left( {n - k}
\right) + \varepsilon }} {{V^{\frac{1}
{{n - k - 1}}} }}} \right)^{ - 1}  - \theta } \right]\\
& = &\frac{{\left( {1 + c\varepsilon _0 } \right)^{\frac{{n - k}}
{{n - k - 1}}} }} {{\left( {1 - c\varepsilon _0 }
\right)^{\frac{n-k+1} {{n - k - 1}}} }}V^{\frac{{n - k}} {{n - k -
1}} - 1 - \frac{2} {{n - k - 1}}} \left( {\frac{{V^{\frac{1} {{n -
k - 1}}} }}
{{\tilde A_0\left( {n - k} \right) + \varepsilon }} - \theta } \right) \hfill \\
  & = &\frac{{\left( {1 + c\varepsilon _0 } \right)^{\frac{{n - k}}
{{n - k - 1}}} }} {{\left( {1 - c\varepsilon _0 }
\right)^{\frac{n-k+1} {{n - k - 1}}} }}V^{\frac{{ - 1}} {{n - k -
1}}} \left( {\frac{{V^{\frac{1} {{n - k - 1}}} }}
{{\tilde A_0\left( {n - k} \right) + \varepsilon }} - \theta } \right) \hfill \\
 & = &\frac{{\left( {1 + c\varepsilon _0 } \right)^{\frac{{n - k}}
{{n - k - 1}}} }} {{\left( {1 - c\varepsilon _0 }
\right)^{\frac{n-k+1} {{n - k - 1}}} }}\left( {\frac{1}
{{\tilde A_0( {n - k}) + \varepsilon }} - \tilde \theta } \right) \hfill \\
\end{eqnarray*}
thus for $\varepsilon_0$  small enough we have
\begin{eqnarray*}
\frac{{\left( {\int_N {\left| {\nabla _{\bar g} u_\alpha }
\right|^2 d\upsilon _{\bar g} }  + \alpha \int_{\partial N}
{u_\alpha ^2 ds_{\bar g} } } \right)\left( {\int_{\partial N}
{\left| {u_\alpha } \right|} ds_{\bar g} } \right)^{\frac{2} {{n -
k - 1}}} }} {{\left( {\int_{\partial N} {u_\alpha ^2 ds_{\bar g} }
} \right)^{\frac{{n - k}} {{n - k - 1}}} }} < \\\left( {1 +
\varepsilon } \right)\left( {\frac{1}
{{\tilde A_0\left( {n - k} \right)}} - \tilde \theta } \right) \hfill \\
\end{eqnarray*}
The latter inequality for $\varepsilon$ small enough yields
\begin{equation}\label{E4.20}
\frac{{\left( {\int_N {\left| {\nabla _{\bar g} u_\alpha }
\right|^2 d\upsilon _{\bar g} }  + \alpha \int_{\partial N}
{u_\alpha ^2 ds_{\bar g} } } \right)\left( {\int_{\partial N}
{\left| {u_\alpha } \right|} ds_{\bar g} } \right)^{\frac{2} {{n -
k - 1}}} }} {{\left( {\int_{\partial N} {u_\alpha ^2 ds_{\bar g} }
} \right)^{\frac{{n - k}} {{n - k - 1}}} }} \leqslant \,\,\frac{1}
{{\tilde A_0\left( {n - k} \right)}}
\end{equation}
Because of (\ref{E4.17}), inequality (\ref{E4.20}) is false and
the proof is complete.\\
\noindent{\bf{Proof of Corollary \ref{T2.3}}}. The proof of the
first part of Corollary \ref{T2.3} arises immediately by Theorem
\ref{T2.2}. (An another proof is presented in \cite{Cot-Lab-Tra}).
For the second part we use Theorem 1 in \cite{Hum2}, since we have
calculated the value of the first Neumann eigenvalue for the
Laplacian on radial functions on the interval $[-1,1]$ .





\begin{thebibliography}{99}





\bibitem{Bak-Cou-Led-Sal}
D, Bakry, T. Coulhon, M. Ledoux and L. Saloff-Coste, \emph{Sobolev
inequalities in disguise},  Indiana University Mathematics
Journal, 44 (1995), 1033-1074.

\bibitem{Bec}
W. Beckner, \emph{Geometric proof of Nash's inequality}, Internat.
Math. Res. Notices (IMRN), (1998), 67-71.


\bibitem{Car-Los}
E. A. Carlen and M. Loss, \emph{Sharp constant in Nash's
inequality},  International Mathematics Research Notices, { 7}
(1993) 213-215.

\bibitem{Cot-Lab3}
A. Cotsiolis  and N. Labropoulos, \emph{Best constants in Sobolev
inequalities on manifolds with boundary in the presence of
symmetries and applications},  Bull. Sci. math. 132 (2008),
562-574.

\bibitem{Cot-Lab-Tra}
A. Cotsiolis, N. Labropoulos and E. Traboulay, \emph{Optimal Nash
Inequalities on the Solid Torus}, USC Journal of Research I
(2008), 47-60.


\bibitem{Dru-Heb-Vau}
O. Druet, E. Hebey and M. Vaugon, \emph{Optimal Nash's
inequalities on Riemannian manifolds: the influence of geometry},
Internat. Math. Res. Notices, no. 14 (1999),  735-779.

\bibitem{Fag1}
Z. Faget, \emph{Best constants in Sobolev inequalities on
Riemannian manifolds in the presence of symmetries}, it Potential
Analysis, { 17} (2002), 105-124.

\bibitem{Fag2}
Z. Faget, \emph{Optimal constants in critical Sobolev inequalities
on Riemannian manifolds in the presence of symmetries}, Annals of
Global Analysis and Geometry, { 24} (2003), 161-200.

\bibitem{Heb}
E. Hebey, "Nonlinear Analysis on Manifolds: Sobolev Spaces and
Inequalities", Courant Institute of Mathematical Sciences,
Lectures Notes in Mathematics, 5, 1999.

\bibitem{Heb-Vau}
E. Hebey and M.Vaugon, \emph{Sobolev spaces in the presence of
symmetries}, J.Math.Pures Appl., { 76} (1997) 859-881.

\bibitem{Hum1}
E. Humbert,  \emph{Best constants in the $L\sp 2$-Nash
inequality}, Proc. Roy. Soc. Edinburgh Sect. A 131  no. { 3},
(2001), 621-646.


\bibitem{Hum2}
E. Humbert, \emph{Optimal trace Nash inequality}, Geom. Funct.
Anal. 11 no. { 4} (2001),  759-772.

\bibitem{Kat}
T. Kato, \emph{The Navier-Stokes equation for an incompressible
fluid in $\mathbb{R}^2$ with a measure as the initial vorticity},
Differential Integral Equations { 7} (1994), 949-966.

\bibitem{Li-Zhu}
Y. Y. Li and M. Zhu,  \emph{Sharp Sobolev inequalities involving
boundary terms},  Geom. Funct. Anal. {\bf 8}  no. 1, (1998),
59--87.

\bibitem{Nas}
J. Nash, \emph{Continuity of solutions of parabolic and eliptic
equations}, American Journal of Mathematics, { 80} (1958),
931-954.
\end{thebibliography}
\end{document}